\documentclass[11pt]{article}
\usepackage{amsmath,amsthm,amsfonts,amssymb}
\usepackage{epsfig}
\usepackage{color}
\usepackage{latexsym}
\usepackage{supertabular}
\usepackage{graphicx}
\usepackage{a4wide}

\numberwithin{equation}{section}

\def\whitebox{{\hbox{\hskip 1pt
 \vrule height 6pt depth 1.5pt
 \lower 1.5pt\vbox to 7.5pt{\hrule width
    3.2pt\vfill\hrule width 3.2pt}%
 \vrule height 6pt depth 1.5pt
 \hskip 1pt } }}
\def\qed{\ifhmode\allowbreak\else\nobreak\fi\hfill\quad\nobreak
     \whitebox\medbreak}

\newcommand{\ignore}[1]{}

\theoremstyle{plain}
\newtheorem{theorem}{Theorem}[section]
\newtheorem{corollary}[theorem]{Corollary}
\newtheorem{lemma}[theorem]{Lemma}

\newtheorem{example}[theorem]{Example}

\newtheorem{remark}[theorem]{Remark}

\def\qed{{\hfill$\square$}}
\def\proof{{\vspace{-0.3cm}\bf Proof: \,}}

\def\Z{{\mathbb Z}}
\def\C{{\mathbb C}}
\def\Q{{\mathbb Q}}

\def\F{{\mathbb F}}
\def\mod{{\mathrm{mod\,\,}}}

\def\Tr{{\mathrm{Tr}}}
\def\Cay{{\mathrm{Cay}}}
\def\ord{{\mathrm{ord}}}

\title{Constructions of Strongly Regular Cayley Graphs and Skew Hadamard Difference Sets from Cyclotomic Classes}
\author{Tao Feng\footnotemark[1]\hspace{10mm} Koji Momihara\footnotemark[2]\hspace{10mm} Qing Xiang\footnotemark[3]}
\date{} 
\begin{document}
\maketitle

\footnotetext[1]{Department of Mathematics, Zhejiang University, Hangzhou 310027, Zhejiang, China;
\\Email: pku.tfeng@yahoo.com.cn}

\footnotetext[2]{
Faculty of Education, Kumamoto University,  
2-40-1 Kurokami, Kumamoto 860-8555, Japan; \\Email: 
momihara@educ.kumamoto-u.ac.jp}

\footnotetext[3]{
Department of Mathematical Science, University of Delaware, Newark, DE 19716, USA; \\Email: xiang@math.udel.edu}
\renewcommand{\thefootnote}{\arabic{footnote}}
%%%%%%%%%%%%%%%%%%%%%%%%%%%%%%%%%%%%%%%%%%%%%%%%%%%%%
%%%%%%%%%%%%%%%%%%%%%%%%%%%%%%%%%%%%%%%%%%%%%%%%%%%%%
\begin{abstract}
In this paper, we give a construction of strongly regular Cayley graphs and a construction of skew Hadamard difference sets. 
Both constructions are based on choosing  cyclotomic classes in finite fields, and they generalize the constructions given by Feng and Xiang \cite{FX111,FX113}. 
Three infinite families of strongly regular graphs with new parameters are obtained. The main tools that we employed are index $2$ Gauss sums, instead of cyclotomic numbers.      
\end{abstract}
\begin{center} 
\end{center}
%%%%%%%%%%%%%%%%%%%%%%%%%%%%%%%%%%%%%%%%%%%%%%%%%%%%%
%%%%%%%%%%%%%%%%%%%%%%%%%%%%%%%%%%%%%%%%%%%%%%%%%%%%%
\section{Introduction}

We assume that the reader is familiar with the basic theory of strongly regular graphs and difference sets. For the theory of strongly regular graphs, 
our main references are  the lecture notes of Brouwer and Haemers \cite{bh} and \cite{cg}. For difference sets, we refer the reader to \cite{Lander} and Chapter 6 of \cite{bjl}. 
We remark that strongly regular graphs are closely related to other combinatorial objects, such as two-weight codes, two-intersection sets in finite geometry, 
and partial difference sets. For these connections, we refer the reader to \cite[p.~132]{bh}, \cite{CK86, M94}. 

Let $\Gamma$ be a (simple, undirected) graph. The adjacency matrix of $\Gamma$ is the $(0,1)$-matrix $A$ with both rows and columns indexed by the vertex set of $\Gamma$, where $A_{xy} = 1$ when there is an edge between $x$ and $y$ in $\Gamma$ and $A_{xy} = 0$ otherwise. A useful way to check whether a graph is strongly regular is by using the eigenvalues of its adjacency matrix. For convenience we call an eigenvalue {\it restricted} if it has an eigenvector perpendicular to the all-ones vector ${\bf 1}$. (For a $k$-regular connected graph, the restricted eigenvalues are the eigenvalues different from $k$.)

\begin{theorem}\label{char}
For a simple graph $\Gamma$ of order $v$, not complete or edgeless, with adjacency matrix $A$, the following are equivalent:
\begin{enumerate}
\item $\Gamma$ is strongly regular with parameters $(v, k, \lambda, \mu)$ for certain integers $k, \lambda, \mu$,
\item $A^2 =(\lambda-\mu)A+(k-\mu) I+\mu J$ for certain real numbers $k,\lambda, \mu$, where $I, J$ are the identity matrix and the all-ones matrix, respectively, 
\item $A$ has precisely two distinct restricted eigenvalues.
\end{enumerate}
\end{theorem}

One of the most effective methods for constructing strongly regular graphs is by the Cayley graph construction. For example, the Paley graph ${\rm P}(q)$ and the Clebsch graph are both Cayley graphs (moreover they are cyclotomic). Let $G$ be an additively written group of order $v$, and let $D$ be a subset of $G$ such that $0\not\in D$ and $-D=D$, where $-D=\{-d\mid d\in D\}$. The {\it Cayley graph on $G$ with connection set $D$}, denoted ${\rm Cay}(G,D)$, is the graph with the elements of $G$ as vertices; two vertices are adjacent if and only if their difference belongs to $D$. In the case when $\Cay(G,D)$ is a strongly regular graph, the connection set $D$ is called a (regular) {\it partial difference set}.  The survey of Ma~\cite{M94} contains much of what is known about partial difference sets and about connections with strongly regular graphs. 

A difference set $D$ in an (additively written) finite group $G$ is called {\it skew Hadamard} if $G$ is the disjoint union of $D$, $-D$, and $\{0\}$. The primary example (and for many years, the only known example in abelian groups) of skew Hadamard difference sets is the classical Paley difference set in $(\F_q,+)$ consisting of the nonzero squares of $\F_q$, where $\F_q$ is the finite field of order $q$, and $q$ is a prime power congruent to 3 modulo 4.  This situation changed dramatically in recent years. Skew Hadamard difference sets are currently under intensive study; see the introduction of \cite{FX113} for a short survey of known constructions of skew Hadamard difference sets and related problems.
 
As we have seen above, in order to obtain strongly regular Cayley graphs, we need to construct regular partial difference sets. 
A classical method for constructing both partial difference sets and difference sets in the additive groups of finite fields is to use cyclotomic classes of finite fields. 
Let $p$ be a prime,  $f$ a positive integer, and let $q=p^f$. Let $N>1$ be an integer such that $N|(q-1)$, and $\gamma$ be a primitive element of $\F_q$. 
Then the cosets $C_i=\gamma^i \langle \gamma^N\rangle$, $0\leq i\leq N-1$, are called the {\it cyclotomic classes of order $N$} of $\F_q$. 
The numbers $|(C_i+1)\cap C_j|$ are called {\it cyclotomic numbers}.  Many authors have studied the problem of determining when a union of some cyclotomic classes forms a (partial) difference set. 
A summary of results in this direction obtained up to 1967 appeared in \cite{S67}. However, all the results in \cite{S67} are based on cyclotomic classes of small orders and 
this method has had only very limited success. In fact, known infinite series of difference sets were obtained only in the case when $N=2$. The situation for partial difference sets is slightly better 
(see \cite{bwx}, \cite[p.~137-138]{bh}). 
One of the reasons why very few difference sets have been discovered by this method is the difficulty of computing cyclotomic numbers of order $N$ when $N$ is large. So far cyclotomic numbers have been evaluated for $N\le 24$ \cite{BEW97} (but note that some of these evaluations are not explicit). For large $N$, probably Van Lint and Schrijver \cite{VLSch} are the first to use cyclotomic classes of order $N$ of finite fields to construct strongly regular graphs, and Baumert, Mills and Ward \cite{BMW82} are the first to use cyclotomic classes of order $N$ of finite fields to construct difference sets. We comment that the difference sets constructed in \cite{BMW82} are also partial difference sets since the finite fields involved have characteristic 2. Both constructions are based on the so-called uniform cyclotomy, which will be defined in Section 2. 

On the other hand, many sporadic examples of strongly regular Cayley graphs have been found using unions of cyclotomic classes of $\F_q$ by computer search. For example, the following are known: 
\begin{itemize}
\item[(i)] (De Lange \cite{DL95}) Let $q=2^{12}$ and  $N=45$. Then,  $\Cay(\F_q, C_0\cup C_5\cup C_{10})$ is a strongly regular graph.

\item[(ii)] (Ikuta and Munemasa \cite{IM10}) Let $q=2^{20}$ and  $N=75$. Then,  
$\Cay(\F_q, C_0\cup C_3\cup C_{6}\cup C_9\cup C_{12})$ is a strongly regular graph.

\item[(iii)] (Ikuta and Munemasa \cite{IM10}) Let $q=2^{21}$ and  $N=49$. Then,  
$\Cay(\F_q, C_0\cup C_1\cup C_{2}\cup C_3\cup C_{4}\cup C_5\cup C_6)$ is a strongly regular graph.

\end{itemize}
Recently, in \cite{FX111}, the first and the third authors extended the above examples to infinite families by using index $2$ Gauss sums over $\F_q$. Below is the main theorem from \cite{FX111}.
\begin{theorem}\label{IntroThm2}
\begin{itemize}
\item[(i)] Let $p_1\equiv 3\,(\mod{4})$ be a prime, $p_1\not=3$, $N=p_1^m$, and let $p$ be a prime 
such that $f:=\ord_N(p)=\phi(N)/2$, where $\phi$ is the Euler totient function. 
Let $q=p^f$ and $D=\bigcup_{i=0}^{p_1^{m-1}-1}C_i\subseteq \F_q$. Assume that 
$1+p_1=4p^h$, where $h$ is the class number of $\Q(\sqrt{-p_1})$. Then $\Cay(\F_q, D)$ is a strongly regular graph.
\item[(ii)]
Let $p_1$ and $p_2$ be primes such that $\{p_1\,(\mod{4}),p_2\,(\mod{4})\}=\{1,3\}$, $N=p_1^mp_2$, and let $p$ be a prime 
such that $\ord_{p_1^m}(p)=\phi(p_1^m)$, $\ord_{p_2}(p)=\phi(p_2^n)$, and $f:=\ord_N(p)=\phi(N)/2$. 
Let $q=p^f$ and $D=\bigcup_{i=0}^{p_1^{m-1}-1}C_{ip_2}\subseteq \F_q$. Assume that 
$p_1=2p^{h/2}+(-1)^{\frac{p_1-1}{2}}b$, $p_2=2p^{h/2}-(-1)^{\frac{p_1-1}{2}}b$, $h$ is even, and $1+p_1p_2=4p^h$, where $b\in \{1,-1\}$ and $h$ is the class number of $\Q(\sqrt{-p_1p_2})$. Then $\Cay(\F_q, D)$ is a strongly regular graph.
\end{itemize}
\end{theorem}
Furthermore, in \cite{FX113}, the following two constructions of skew Hadamard  difference sets and Paley type partial difference sets were given. 
(A partial difference set $D$ in a group $G$ is said to be of {\it Paley type} if the parameters of the corresponding strongly regular Cayley graph are $(v, \frac{v-1}{2}, \frac{v-5}{4}, \frac{v-1}{4})$.)

\begin{theorem}\label{IntroThm3}
\begin{itemize}
\item[(i)] Let $p_1\equiv 7\,(\mod{8})$ be a prime, $N=2p_1^m$, and 
let $p$ be a prime such that $f:=\ord_N(p)=\phi(N)/2$. Let $s$ be an odd integer, $I$ any subset of $\Z_N$ such that 
$\{i\,(\mod{p_1^m})\,|\,i\in I\}=\Z_{p_1^m}$, and  let $D=\bigcup_{i\in I}C_i\subseteq \F_{p^{fs}}$. Then, $D$ is a skew Hadamard 
difference set if $p\equiv 3\,(\mod{4})$ and $D$ is a 
Paley type partial difference set if $p\equiv 1\,(\mod{4})$. 
\item[(ii)] Let $p_1\equiv 3\,(\mod{8})$ be a prime, $p_1\neq 3$, $N=2p_1$, and 
let $p\equiv 3\,(\mod{4})$ be a prime such that $f:=\ord_N(p)=\phi(N)/2$. 
Let $q=p^f$, $I=\langle p\rangle\cup 2\langle p\rangle\cup\{0\}$, and let $D=\bigcup_{i\in I}C_i\subseteq \F_q$. Assume that $1+p_1=4p^h$, where $h$ is the class number of $\Q(\sqrt{-p_1})$. Then, $D$ is a skew Hadamard 
difference set  in the additive group of $\F_q$. 
\end{itemize}
\end{theorem}
Note that in Theorem~\ref{IntroThm3} (ii), we need to choose a suitable primitive element $\gamma$ of $\F_q$. For details,  see \cite{FX113}. To extend the construction of Theorem~\ref{IntroThm3} (ii) to the general case $N=2p_1^m$ was left as an open problem in \cite{FX113}. 

The purpose of this paper is to generalize the constructions of strongly regular Cayley graphs in Theorem~\ref{IntroThm2} (ii) to the case where $N=p_1^mp_2^n$ and 
of skew Hadamard difference sets in Theorem~\ref{IntroThm3} (ii) to 
the case where $N=2p_1^m$. Three infinite families of strongly regular graphs with new parameters are obtained (see Table 2 in Section 3). 
An infinite series of skew Hadamard difference sets in $(\F_q,+)$, where $q=3^{53\cdot 107^{m-1}}$, is also obtained. Implications of these results on association schemes will be discussed in Section 4.
%%%%%%%%%%%%%%%%%%%%%%%%%%%%%%%%%%%%%%%%%%%%%%%%%%%%%%%%%%%%%%%%%%%%%%%%%%%%%%%%%%%%%%%
%%%%%%%%%%%%%%%%%%%%%%%%%%%%%%%%%%%%%%%%%%%%%%%%%%%%%%%%%%%%%%%%%%%%%%%%%%%%%%%%%%%%%%%
\section{Index 2 Gauss sums}
Let $p$ be a prime, $f$ a positive integer, and $q=p^f$. The canonical additive character $\psi$ of $\F_q$ is defined by 
$$\psi\colon\F_q\to \C^{*},\qquad\psi(x)=\zeta_p^{\Tr _{q/p}(x)},$$
where $\zeta_p={\rm exp}(\frac {2\pi i}{p})$ and $\Tr _{q/p}$ is the trace from $\F_q$ to $\F_p$. For a multiplicative character 
$\chi$ of $\F_q$, we define the {\it Gauss sum} 
\[
G(\chi)=\sum_{x\in \F_q^\ast}\chi(x)\psi(x). 
\] 
%In this paper, for convenience, we extend the domain of $\chi$ to all elements of $\F_q$ by defining %$\chi(0)=1$ or $0$ depending on whether 
%$\chi$ is trivial or not. 

Below are a few basic properties of Gauss sums \cite{LN97}: 
\begin{itemize}
\item[(i)] $G(\chi)\overline{G(\chi)}=q$ if $\chi$ is nontrivial;
\item[(ii)] $G(\chi^p)=G(\chi)$, where $p$ is the characteristic of $\F_q$; 
\item[(iii)] $G(\chi^{-1})=\chi(-1)\overline{G(\chi)}$;
\item[(iv)] $G(\chi)=-1$ if $\chi$ is trivial. 
\end{itemize}

In general, the explicit evaluation of Gauss sums is a very difficult problem. There are only a few cases where the Gauss sums have been evaluated. 
The simplest case is the so-called {\it semi-primitive case} (also 
referred to as {\it uniform cyclotomy} or {\it pure Gauss sum}), where there 
exists an integer $j$ such that $p^j\equiv -1\,(\mod{N})$, here $N$ is the order of 
the multiplicative character $\chi$ involved. See \cite{BEW97,BMW82,CK86} for the explicit evaluation in this case. 

The next interesting case is the index $2$ case where the subgroup $\langle p\rangle$ generated by $p\in \Z_{N}^\ast$ has index $2$ in $\Z_{N}^\ast$ and $-1\not\in \langle p\rangle $. In this case, 
it is known that $N$ can have at most two odd prime divisors. 
Many authors have investigated this case, see e.g., \cite{BM73,L97,M98,MV03,YX10,YX11}. In particular, a complete solution to the problem of evaluating Gauss sums in this case was recently given in \cite{YX10}. The following are the results on evaluation of Gauss sums which 
we will need in the next section. 
\begin{theorem}\label{Sec2Thm1}(\cite{YX10},  Case B1; Theorem~4.10)
Let $N=p_1^mp_2^n$, where $m$ and $n$ are positive integers, $p_1$ and $p_2$ are  primes such that $p_1\equiv 1\,(\mod{4})$ and $p_2\equiv 3\,(\mod{4})$. Assume that $p$ is a prime such that 
$\ord_{p_1^m}(p)=\phi(p_1^m)$, $\ord_{p_2^n}(p)=\phi(p_2^n)$, and $[\Z_N^\ast:\langle p \rangle]=2$. 
Let $f=\phi(N)/2$, $q=p^f$, and $\chi$ be a multiplicative character of order $N$ of $\F_q$. Then, for 
$0\le s\le m-1$ and $0\le t\le n-1$, we have
\begin{eqnarray*}
G(\chi^{p_1^sp_2^t})&=&p^{\frac{f-hp_1^sp_2^t}{2}}
\left(\frac{b+c\sqrt{-p_1p_2}}{2}\right)^{p_1^sp_2^t};\\
G(\chi^{p_1^mp_2^t})&=&-p^{\frac{f}{2}};\\
G(\chi^{p_1^sp_2^n})&=&p^{\frac{f}{2}}, 
\end{eqnarray*}
where $h$ is the class number of $\Q(\sqrt{-p_1p_2})$, and $b$ and $c$ are integers 
determined by $b,c\not\equiv 0\,(\mod{p})$, $4p^{h}=b^2+p_1p_2c^2$,  and $bp^{\frac{f-h}{2}}\equiv 2\,(\mod{p_1p_2})$. 
\end{theorem}

\begin{theorem}\label{Sec2Thm2}(\cite{YX10}, Case D; Theorem~4.12)
Let $N=2p_1^m$, where $p_1>3$ is a prime such that $p_1\equiv 3\,(\mod{4})$ and $m$ is a positive integer. Assume that $p$ is a 
prime such that  $[\Z_N^\ast:\langle p \rangle]=2$. 
Let $f=\phi(N)/2$, $q=p^f$, and $\chi$ be a multiplicative character of order $N$ of $\F_q$. Then, for 
$0\le t\le m-1$, we have
\begin{eqnarray*}
G(\chi^{p_1^t})&=&\left\{
\begin{array}{ll}
(-1)^{\frac{p-1}{2}(m-1)}p^{\frac{f-1}{2}-hp_1^t}
\sqrt{p^\ast}\left(\frac{b+c\sqrt{-p_1}}{2}\right)^{2p_1^t},&  \mbox{if $p_1\equiv 3\,(\mod{8})$,}\\
(-1)^{\frac{p-1}{2}m}p^{\frac{f-1}{2}}
\sqrt{p^\ast},&  \mbox{if $p_1\equiv 7\,(\mod{8})$;}
 \end{array}
\right.\\
G(\chi^{2p_1^t})&=&p^{\frac{f-p_1^t h}{2}}\left(\frac{b+c\sqrt{-p_1}}{2}\right)^{p_1^t};\\
G(\chi^{p_1^m})&=&(-1)^{\frac{p-1}{2}\frac{f-1}{2}}p^{\frac{f-1}{2}}\sqrt{p^\ast}, 
\end{eqnarray*}
where $p^\ast=(-1)^{\frac{p-1}{2}}p$, $h$ is the class number of $\Q(\sqrt{-p_1})$, and $b$ and $c$ are integers 
determined by $4p^{h}=b^2+p_1c^2$ and $bp^{\frac{f-h}{2}}\equiv -2\,(\mod{p_1})$. 
\end{theorem}

Note that Theorem~\ref{Sec2Thm2} above is Theorem~4.12 in \cite{YX10}, whose statement contains several misprints. We corrected those misprints in the above statement.
 
%%%%%%%%%%%%%%%%%%%%%%%%%%%%%%%%%%%%%%%%%%%%%%%%%%%%%%%%%%%%%%%%%%
%%%%%%%%%%%%%%%%%%%%%%%%%%%%%%%%%%%%%%%%%%%%%%%%%%%%%%%%%%%%%%%%%%
\section{Constructions of strongly regular Cayley graphs and skew Hadamard difference sets}
We first recall the following well-known lemma in the theory of difference sets 
(see e.g., \cite{M94,T65}). 
\begin{lemma}\label{Sec3Le1}
Let $(G, +)$ be an abelian group of odd order $v$, $D$ a subset of $G$ 
of size $\frac{v-1}{2}$.  Assume that $D\cap -D=\emptyset$ and 
$0\not \in D$. Then, $D$ is a skew Hadamard difference set 
in $G$ if and only if 
\[
\chi(D)=\frac{-1\pm \sqrt{-v}}{2}
\] 
for all nontrivial characters $\chi$ of $G$.  On the other hand, assume that $0\not\in D$ and $-D=D$. Then $D$ is a Paley type partial difference set in $G$ if and only if 
\[
\chi(D)=\frac{-1\pm \sqrt{v}}{2}
\] 
for all nontrivial characters $\chi$ of $G$. 
\end{lemma}

Let $q=p^f$, where $p$ is a prime and $f$ a positive integer, and let 
$C_i=\gamma^i \langle \gamma^N\rangle$, $0\le i\le N-1$, be 
the cyclotomic classes of order $N$ of $\F_q$, where $\gamma$ is a fixed primitive element of $\F_q$. From now on, we will assume that 
$D$ is a union of cyclotomic classes of order $N$ of $\F_q$. In order to check whether a candidate subset, $D=\bigcup_{i\in I}C_i$, is a partial difference set or a skew Hadamard difference set in $(\F_q,+)$, we will compute the sums $\psi(aD):=\sum_{x\in D}\psi(ax)$ for all $a\in \F_q^\ast$,  where $\psi$ is the canonical additive character of $\F_q$. Note that the sum $\psi(aD)$ can be expressed as a linear combination of 
Gauss sums using the orthogonality of characters: 
\begin{eqnarray*}
\psi(aD)&=&\frac{1}{N}\sum_{i\in I}\sum_{x\in \F_q^\ast}\psi(a\gamma^i x^N)\\
&=&\frac{1}{N}\sum_{i\in I}\sum_{x\in \F_q^\ast}\frac{1}{q-1}
\sum_{y\in \F_q^\ast}\psi(y)
\sum_{\chi\in \widehat{\F_q^\ast}}\chi(a\gamma^ix^N)\overline{\chi(y)}\\
&=&\frac{1}{(q-1)N}\sum_{i\in I}\sum_{x\in \F_q^\ast}
\sum_{\chi\in \widehat{\F_q^\ast}}G(\chi^{-1})
\chi(a\gamma^ix^N)\\
&=&\frac{1}{(q-1)N}\sum_{i\in I}
\sum_{\chi\in \widehat{\F_q^\ast}}G(\chi^{-1})
\chi(a\gamma^i)\sum_{x\in \F_q^\ast}\chi(x^N)\\
&=&\frac{1}{N}
\sum_{\chi\in C_0^{\perp}}G(\chi^{-1})
\sum_{i\in I}\chi(a\gamma^i ), 
\end{eqnarray*}
where $\widehat{\F_q^\ast}$ is the group of multiplicative characters of  
$\F_q^\ast$ and 
$C_0^{\perp}$ is the subgroup of $\widehat{\F_q^\ast}$
consisting of all $\chi$ which are trivial on $C_0$. 
%%%%%%%%%%%%%%%%%%%%%%%%%%%%%%%%%%%%%%%%%%%%%%%%%%%%%%%%%%%%%%%%%%%%%%%%%%%%%%
%%%%%%%%%%%%%%%%%%%%%%%%%%%%%%%%%%%%%%%%%%%%%%%%%%%%%%%%%%%%%%%%%%%%%%%%%%%%%%
\subsection{Strongly regular graphs from unions of cyclotomic classes of order $N=p_1^mp_2^n$}\label{Subsec1Tit}
In this subsection, we assume 
that $N=p_1^mp_2^n$, where $m, n$ are positive integers, $p_1$ and $p_2$ are primes such that $p_1\equiv 1\,(\mod{4})$ and $p_2\equiv 3\,(\mod{4})$. 
Furthermore, we assume that $p$ is a prime such that  $\ord_{p_1^m}(p)=\phi(p_1^m)$,  $\ord_{p_2^n}(p)=\phi(p_2^n)$, and $\ord_N(p)=\phi(N)/2$. Let $q=p^f$ 
and $C_i=\gamma^i \langle \gamma^N\rangle$, $0\leq i\leq N-1$, where $f=\ord_N(p)$ and  $\gamma$ is a 
fixed primitive element of $\F_q$.  
Define 
\[
D=\bigcup_{i=0}^{p_1^{m-1}-1}\bigcup_{j=0}^{p_2^{n-1}-1}C_{p_2^n i+p_1^{m}j}. 
\] 
It is clear that $D=-D$.
\begin{theorem}\label{Subsec1Thm1}
The size of the set $\{\psi(\gamma^a D)\,|\,a=0,1,\ldots,q-2\}$ is at most five.  
\end{theorem}
\proof 
Let $\chi_{e}$ denote the multiplicative character of 
order $e$ of $\F_q$ such that 
$\chi_e(\gamma)=\zeta_e$, where $\zeta_e:=\exp{(\frac{2\pi i}{e})}$. 
Then $\chi_{e}^{d}=\chi_{\frac{e}{d}}$ for any divisor $d$ of $e$.  Note that since $D$ is a union of cyclotomic classes of order $N$, we have $\{\psi(\gamma^a D)\,|\,a=0,1,\ldots,q-2\}=\{\psi(\gamma^a D)\,|\,a=0,1,\ldots, N-1\}$.

To prove the theorem,  
it is sufficient to evaluate the sums 
\[
T_a:=N\cdot \psi(\gamma^a D)=\sum_{\ell=0}^{p_1^m p_2^n-1}G(\chi_{p_{1}^m p_2^n}^{-\ell}) \sum_{i=0}^{p_1^{m-1}-1}\sum_{j=0}^{p_2^{n-1}-1}
\chi_{p_{1}^mp_2^n}^{\ell}(\gamma^{a+p_2^ni+p_1^{m}j}), 
\]
 where $a=0,1,\ldots, N-1$.

For $\ell=0$, by noting that $G(\chi_{p_{1}^mp_2^n}^{0})=-1$, we have 
\[
G(\chi_{p_{1}^mp_2^n}^{0}) \sum_{i=0}^{p_1^{m-1}-1}\sum_{j=0}^{p_2^{n-1}-1}
\chi_{p_{1}^mp_2^n}^{\ell}(\gamma^{a+p_2^ni+p_1^{m}j})=-p_1^{m-1}p_2^{n-1}. 
\]

For $\ell=p_1 h$ but  $h\not\equiv 0\,(\mod{p_1^{m-1}})$,  
we have 
\[
\sum_{i=0}^{p_1^{m-1}-1}
\chi_{p_{1}^mp_2^n}^{p_1h}(\gamma^{p_2^ni})
=\sum_{i=0}^{p_1^{m-1}-1}
\chi_{p_{1}^{m-1}}^{h}(\gamma^{i})
=0. 
\]

For $\ell=p_2 h$ but $h\not\equiv 0\,(\mod{p_2^{n-1}})$,  
we have 
\[
\sum_{j=0}^{p_2^{n-1}-1}
\chi_{p_{1}^mp_2^n}^{p_2h}(\gamma^{p_1^{m}j})
=\sum_{j=0}^{p_2^{n-1}-1}
\chi_{p_2^{n-1}}^{h}(\gamma^{j})=0. 
\]

Note that for each $a\in \{0,1,\ldots,N-1\}$, there is a unique 
$i\in \{0,1,\ldots,p_1^{m-1}-1\}$ such that 
$p_1^{m-1}\,|\,a+p_2^n i$; we write $a+p_2^n i=p_1^{m-1}i_a$. 
Define $\delta_{i_a}=1$ or $0$ depending on whether $i_a\equiv 0\,(\mod{p_1})$ or 
not. Similarly, for each $a\in \{0,1,\ldots,N-1\}$, there is a unique 
$j\in \{0,1,\ldots,p_2^{n-1}-1\}$ such that 
$p_2^{n-1}\,|\,a+p_1^m j$; we write $a+p_1^m j=p_2^{n-1}j_a$. 
Define $\delta_{j_a}=1$ or $0$ depending on whether $j_a\equiv 0\,(\mod{p_2})$ or not. 

For $\ell=p_1^m h$ but  $h\not\equiv 0\,(\mod{p_2})$, since 
$G(\chi_{p_{2}^n}^{-h})=-p^{\frac{f}{2}}$ by Theorem~\ref{Sec2Thm1}, 
we have 
\begin{eqnarray*}
& &\sum_{h:\gcd{(h,p_2)}=1}G(\chi_{p_{2}^n}^{-h}) \sum_{i=0}^{p_1^{m-1}-1}\sum_{j=0}^{p_2^{n-1}-1}
\chi_{p_2^n}^{h}(\gamma^{a+p_2^ni+p_1^{m}j})\\
&=&-p_1^{m-1}p^{\frac{f}{2}}\sum_{x\in \Z_{p_2}^\ast}
\sum_{y=0}^{p_2^{n-1}-1} \sum_{j=0}^{p_2^{n-1}-1}
\chi_{p_2^n}^{x+p_2y}(\gamma^{a+p_1^{m}j})\\
&=&-p_1^{m-1}p_2^{n-1}p^{\frac{f}{2}}\sum_{x\in \Z_{p_2}^\ast}
\chi_{p_2^n}^{x}(\gamma^{j_a})\\
&=&-p_1^{m-1}p_2^{n-1}p^{\frac{f}{2}}(p_2\delta_{j_a}-1). 
\end{eqnarray*} 
Similarly, for $\ell=p_2^n h$ but $h\not\equiv 0\,(\mod{p_1})$,  
since
$G(\chi_{p_{1}^m}^{-h})=p^{\frac{f}{2}}$ by Theorem~\ref{Sec2Thm1},
we have 
\[
\sum_{h:\gcd{(h,p_1)}=1}G(\chi_{p_{1}^m}^{-h}) \sum_{i=0}^{p_1^{m-1}-1}\sum_{j=0}^{p_2^{n-1}-1}
\chi_{p_1^m}^{h}(\gamma^{a+p_2^ni+p_1^{m}j})
=p_1^{m-1}p_2^{n-1}p^{\frac{f}{2}}(p_1\delta_{i_a}-1). 
\]

For the remaining cases, we consider the sum  
\begin{equation}\label{Subsec1Eq1}
\sum_{\ell:\gcd{(\ell,p_1p_2)}=1}
G(\chi_{p_{1}^mp_2^n}^{-\ell}) \sum_{i=0}^{p_1^{m-1}-1}\sum_{j=0}^{p_2^{n-1}-1}
\chi_{p_{1}^mp_2^n}^{\ell}(\gamma^{a+p_2^ni+p_1^{m}j}). 
\end{equation}
Note that any multiplicative character of order $p_1^mp_2^n$ can be written as 
$\chi_{p_1^m}^u\chi_{p_2^n}^v$ for some $u\in \Z_{p_1^m}^\ast$ and 
$v\in \Z_{p_2^n}^\ast$.  
By Theorem~\ref{Sec2Thm1},  we have 
\[
G(\chi_{p_1^m}^{-u}\chi_{p_2^n}^{-v})=p^{\frac{f-h}{2}}
\frac{b+c\eta_1(u)\eta_2(v)\sqrt{-p_1p_2}}{2}, 
\]
where $b,c\not \equiv 0\,(\mod{p})$, 
$b^2+p_1p_2c^2=4p^h$, $bp^{\frac{f-h}{2}}\equiv 2\,(\mod{p_1p_2})$, 
and $\eta_1$ and $\eta_2$ are the quadratic characters of $\F_{p_1}^\ast$ and $\F_{p_2}^\ast$, respectively. Then, 
the sum (\ref{Subsec1Eq1}) is rewritten as 
\begin{eqnarray} \nonumber 
& &p^{\frac{f-h}{2}}
\sum_{u\in \Z_{p_1^m}^\ast}\sum_{v\in \Z_{p_2^n}^\ast}
\frac{b+c\eta_1(u)\eta_2(v)\sqrt{-p_1p_2}}{2}
\sum_{i=0}^{p_1^{m-1}-1}\sum_{j=0}^{p_2^{n-1}-1}
\chi_{p_{1}^m}^{u}(\gamma^{a+p_2^ni+p_1^{m}j})\chi_{p_2^n}^{v}(\gamma^{a+p_2^ni+p_1^{m}j})\\ 
&=&
\frac{p^{\frac{f-h}{2}}b}{2}\left(\sum_{u\in \Z_{p_1^m}^\ast}
\sum_{i=0}^{p_1^{m-1}-1}\chi_{p_{1}^m}^{u}(\gamma^{a+p_2^ni})\right)
\left(
\sum_{v\in \Z_{p_2^n}^\ast}
\sum_{j=0}^{p_2^{n-1}-1}\chi_{p_2^n}^{v}(\gamma^{a+p_1^{m}j})
\right) \label{Subsec1Eq2}\\ 
& &\hspace{-6mm}+
\frac{p^{\frac{f-h}{2}}c\sqrt{-p_1p_2}}{2}\left(\sum_{u\in \Z_{p_1^m}^\ast}\eta_1(u)
\sum_{i=0}^{p_1^{m-1}-1}\chi_{p_{1}^m}^{u}(\gamma^{a+p_2^ni})\right)
\left(\sum_{v\in \Z_{p_2^n}^\ast}\eta_2(v)
\sum_{j=0}^{p_2^{n-1}-1}
\chi_{p_2^n}^{v}(\gamma^{a+p_1^{m}j})\right) \label{Subsec1Eq3}
\end{eqnarray}

For (\ref{Subsec1Eq2}), we have 
\begin{eqnarray*}
& &\frac{p^{\frac{f-h}{2}}b}{2}\left(\sum_{u\in \Z_{p_1^m}^\ast}
\sum_{i=0}^{p_1^{m-1}-1}\chi_{p_{1}^m}^{u}(\gamma^{a+p_2^ni})\right)
\left(
\sum_{v\in \Z_{p_2^n}^\ast}
\sum_{j=0}^{p_2^{n-1}-1}\chi_{p_2^n}^{v}(\gamma^{a+p_1^{m}j})
\right)\\
&=&\frac{p^{\frac{f-h}{2}}b}{2}\left(\sum_{x\in \Z_{p_1}^\ast}
\sum_{y=0}^{p_1^{m-1}-1}
\sum_{i=0}^{p_1^{m-1}-1}\chi_{p_{1}^m}^{x+p_1y}(\gamma^{a+p_2^ni})\right)
\left(
\sum_{x'\in \Z_{p_2}^\ast}\sum_{y'=0}^{p_2^{n-1}-1}
\sum_{j=0}^{p_2^{n-1}-1}\chi_{p_2^n}^{x'+p_2y'}(\gamma^{a+p_1^{m}j})
\right)\\
&=&\frac{p^{\frac{f-h}{2}}b}{2}\left(p_1^{m-1}\sum_{x\in \Z_{p_1}^\ast}
\chi_{p_{1}}^{x}(\gamma^{i_a})\right)
\left(p_2^{n-1}
\sum_{x'\in \Z_{p_2}^\ast}
\chi_{p_2}^{x'}(\gamma^{j_a})
\right)\\
&=&\frac{p^{\frac{f-h}{2}}b}{2}p_1^{m-1}p_2^{n-1}
(p_1\delta_{i_a}-1)
(p_2\delta_{j_a}-1). 
\end{eqnarray*}

Let $G(\eta_i)$, $i=1,2$, be the quadratic Gauss sums of $\F_{p_i}$, respectively. It is well known that $G(\eta_i)=\sqrt{(-1)^{(p_i-1)/2}p_i}$ (see \cite{LN97}). 
Then, for (\ref{Subsec1Eq3}), we have 
\begin{eqnarray*}
& &\frac{p^{\frac{f-h}{2}}c\sqrt{-p_1p_2}}{2}\left(\sum_{x\in \Z_{p_1}^\ast}
\sum_{y=0}^{p_1^{m-1}-1}\eta_1(x)
\sum_{i=0}^{p_1^{m-1}-1}\chi_{p_{1}^m}^{x+p_1y}(\gamma^{a+p_2^ni})\right)
\left(\sum_{x'\in \Z_{p_2}^\ast}\sum_{y'=0}^{p_2^{n-1}-1}
\eta_2(x')
\sum_{j=0}^{p_2^{n-1}-1}
\chi_{p_2^n}^{x'+p_2y'}(\gamma^{a+p_1^{m}j})\right)\\
&=& \frac{p^{\frac{f-h}{2}}c\sqrt{-p_1p_2}}{2}\left(p_{1}^{m-1}
\sum_{x\in \Z_{p_1}^\ast}
\eta_1(x)\chi_{p_{1}}^{x}(\gamma^{i_a})\right)
\left(p_2^{n-1}\sum_{x'\in \Z_{p_2}^\ast}
\eta_2(x')
\chi_{p_2}^{x'}(\gamma^{j_a})\right)\\
&=& \frac{p^{\frac{f-h}{2}}c\sqrt{-p_1p_2}}{2}p_1^{m-1}p_2^{n-1}
\eta_1(i_a)\eta_2(j_a)G(\eta_1)G(\eta_2)\\
&=& \frac{p^{\frac{f-h}{2}}c\sqrt{-p_1p_2}}{2}p_1^{m-1}p_2^{n-1}
\eta_1(i_a)\eta_2(j_a)\sqrt{(-1)^{\frac{p_1-1}{2}+\frac{p_2-1}{2}}p_1p_2}\\
&=& -\frac{p^{\frac{f-h}{2}}c}{2}p_1^{m}p_2^{n}
\eta_1(i_a)\eta_2(j_a). 
\end{eqnarray*}

Thus, we obtain 
\begin{eqnarray*}
T_a+p_1^{m-1}p_2^{n-1}&=&p_1^{m-1}p_2^{n-1}p^{\frac{f}{2}}
(-p_2\delta_{j_a}+p_1\delta_{i_a})+\frac{p^{\frac{f-h}{2}}b}{2}p_1^{m-1}p_2^{n-1}
(p_1\delta_{i_a}-1)
(p_2\delta_{j_a}-1)\\
& &\hspace{0.5cm}-\frac{p^{\frac{f-h}{2}}c}{2}p_1^{m}p_2^{n}
\eta_1(i_a)\eta_2(j_a). 
\end{eqnarray*}

Now, we compute 
$S_a:=(T_a+p_1^{m-1}p_2^{n-1})/p_1^{m-1}p_2^{n-1}p^{\frac{f}{2}}$ by considering the following four cases:
\begin{itemize}
\item[(i)] If $\delta_{j_a}=\delta_{i_a}=0$, we have 
$
S_a=\frac{p^{\frac{-h}{2}}b}{2}
\pm \frac{p^{\frac{-h}{2}}c}{2}p_1p_2$.
\item[(ii)] If $\delta_{j_a}=1,\delta_{i_a}=0$, we have 
$
S_a=-p_2-\frac{p^{\frac{-h}{2}}b}{2}(p_2-1).
$
\item[(iii)] if $\delta_{j_a}=0,\delta_{i_a}=1$, we have 
$
S_a=p_1-\frac{p^{\frac{-h}{2}}b}{2}(p_1-1). 
$
\item[(iv)] if $\delta_{j_a}=\delta_{i_a}=1$, we have 
$
S_a=-p_2+p_1+\frac{p^{\frac{-h}{2}}b}{2}(p_1-1)(p_2-1). 
$
\end{itemize}
The proof is now complete. 
\qed

\begin{corollary}\label{Subsec1Cor1}
If $b,c\in \{1,-1\}$, 
$h$ is even, $p_1=2p^{h/2}+b$, and 
$p_2=2p^{h/2}-b$, then $\Cay(\F_q,D)$ is a strongly regular graph. 
\end{corollary}
\proof Since $-D=D$ and $0\not\in D$, the Cayley graph ${\rm Cay}(\F_q, D)$ is undirected and without loops. It is also regular of valency $|D|$. 
The restricted eigenvalues of this Cayley graph, as explained in \cite[p.~134]{bh}, are $\psi(\gamma^a D),$ where $a=0,1,\ldots, q-2$. By Theorem~\ref{char}, it suffices to show that the set $\{\psi(\gamma^a D)\,|\,a=0,1,\ldots,q-2\}$ has precisely two elements.  We substitute $p_1=2p^{h/2}+b$,  
$p_2=2p^{h/2}-b$, and $b,c\in \{1,-1\}$ into the expressions for $S_a$ in the proof of Theorem~\ref{Subsec1Thm1}, and find that $S_a$ indeed take only two distinct values. 
This proves that $\Cay(\F_q,D)$ is a strongly regular graph. In particular, the two restricted eigenvalues $r$ and $s$ ($r>s$) are given by 
$r=\frac{2p^{\frac{f+h}{2}}-1}{p_1p_2}$ and 
$s=\frac{-2p^{\frac{f+h}{2}}+p^{\frac{f-h}{2}}-1}{p_1p_2}$,  or 
$s=\frac{-2p^{\frac{f+h}{2}}-1}{p_1p_2}$ and 
$r=\frac{2p^{\frac{f+h}{2}}-p^{\frac{f-h}{2}}-1}{p_1p_2}$ depending on whether 
$b=1$ or $b=-1$. Furthermore, the parameters $k,\lambda,$ and $\mu$ of 
the strongly regular graph are given by $k=\frac{p^{f}-1}{p_1p_2}$, 
$\lambda=s+r+k+sr$, and $\mu=k+sr$. 
\qed

\begin{remark}
One can show that the assumptions on $p_1,p_2,b,c,$ and $h$ are also 
necessary for $\Cay(\F_q,D)$ to be strongly regular by a similar proof to that of Corollary 5.2 in \cite{FX111}.
\end{remark}
The construction of strongly regular Cayley graphs given in this subsection is a generalization of 
Theorem~\ref{IntroThm2} (ii) \cite{FX111}. In \cite{FX111}, the six infinite series of strongly regular graphs in Table~\ref{Subsec1Tab1} below were 
obtained. Note that the case when $m=2$ of the $1$st series of Table~\ref{Subsec1Tab1} is the example found by 
De Lange~\cite{DL95} and 
the case when $m=2$ of the $2$nd series of Table~\ref{Subsec1Tab1} is the example found by Ikuta and Munemasa~\cite{IM10}. 
These six infinite series are combined and generalized to three infinite families of strongly regular graphs in Table 2.
\begin{table}[h]
\begin{center}
\caption{\label{Subsec1Tab1}
Some strongly regular graphs obtained in \cite{FX111}. The parameters $r,s$ are the two nontrivial eigenvalues 
of $\Cay(G,D)$, i.e., the two values in $\{\psi(\gamma^a D)\,|\,a=0,1,\ldots,q-2\}$. 
The parameters $\lambda$ and $\mu$ of the strongly regular graphs can be computed by 
$\lambda=s+r+sr+k$ and $\mu=k+sr$. }
$$
\begin{array}{|c||c|c|c|c|c|c|c|c|c|}
\hline
\mbox{No.}&\mbox{$p$} & \mbox{$N$}&\mbox{$h$}& \mbox{$b$} & \mbox{$v$}&\mbox{$k$}
&\mbox{$r$}&\mbox{$s$}\\
\hline \hline 
\mbox{1}&\mbox{$2$} & \mbox{$3^m\cdot 5$}&\mbox{$2$}& \mbox{$1$} & \mbox{$2^{4\cdot 3^{m-1}}$} & \mbox{$\frac{2^{4\cdot 3^{m-1}}-1}{15}$}
&\mbox{$\frac{8\cdot 2^{2\cdot 3^{m-1}-1}-1}{15}$}&\mbox{$\frac{-7\cdot 2^{2\cdot 3^{m-1}-1}-1}{15}$}\\
\hline
\mbox{2}&\mbox{$2$} & \mbox{$5^m\cdot 3$}&\mbox{$2$}& \mbox{$1$} & \mbox{$2^{4\cdot 5^{m-1}}$} & \mbox{$\frac{2^{4\cdot 5^{m-1}}-1}{15}$}
&\mbox{$\frac{8\cdot 2^{2\cdot 5^{m-1}-1}-1}{15}$}&\mbox{$\frac{-7\cdot 2^{2\cdot 5^{m-1}-1}-1}{15}$}\\
\hline
\mbox{3}&\mbox{$3$} & \mbox{$5^m\cdot 7$}&\mbox{$2$}& \mbox{$-1$} & \mbox{$3^{12\cdot 5^{m-1}}$} & \mbox{$\frac{3^{12\cdot 5^{m-1}}-1}{35}$}
&\mbox{$\frac{17\cdot 3^{6\cdot 5^{m-1}-1}-1}{35}$}&\mbox{$\frac{-18\cdot 3^{6\cdot 5^{m-1}-1}-1}{35}$}\\
\hline
\mbox{4}&\mbox{$3$} & \mbox{$7^m\cdot 5$}&\mbox{$2$}&\mbox{$-1$}& \mbox{$3^{12\cdot 7^{m-1}}$} & \mbox{$\frac{3^{12\cdot 7^{m-1}}-1}{35}$}
&\mbox{$\frac{17\cdot 3^{6\cdot 7^{m-1}-1}-1}{35}$}&\mbox{$\frac{-18\cdot 3^{6\cdot 7^{m-1}-1}-1}{35}$}\\
\hline
\mbox{5}&\mbox{$3$}&\mbox{$17^m\cdot 19$}&\mbox{$4$}&\mbox{$-1$}&\mbox{$3^{144\cdot 17^{m-1}}$} & \mbox{$\frac{3^{144\cdot 17^{m-1}}-1}{323}$}&\mbox{$\frac{161\cdot 3^{72\cdot 17^{m-1}-2}-1}{323}$}&\mbox{$\frac{-162\cdot 3^{72\cdot 17^{m-1}-2}-1}{323}$}\\
\hline
\mbox{6}&\mbox{$3$}&\mbox{$19^m\cdot 17$}&\mbox{$4$}&\mbox{$-1$}&\mbox{$3^{144\cdot 19^{m-1}}$} & \mbox{$\frac{3^{144\cdot 19^{m-1}}-1}{323}$}&\mbox{$\frac{161\cdot 3^{72\cdot 17^{m-1}-2}-1}{323}$}&\mbox{$\frac{-162\cdot 3^{72\cdot 17^{m-1}-2}-1}{323}$}\\
\hline
\end{array}
$$
%\hspace{0.8cm}
\end{center}
\end{table}
\begin{example}\label{Subsec1Ex1}
Table~\ref{Subsec1Tab2} gives generalizations of strongly 
 regular graphs in Table~\ref{Subsec1Tab1}. 
Here, 
the parameters $p,N,h,b$ satisfy the conditions of Corollary~\ref{Subsec1Cor1}, i.e., $p$ is a prime such that 
$[\Z_{N}^\ast:\langle p\rangle]=2$, $b,c\in \{1,-1\}$, 
$p_1=2p^{h/2}+b$, 
$p_2=2p^{h/2}-b$, $h\equiv 0\,(\mod{2})$, and 
$bp^{\frac{f-h}{2}}\equiv 2\,(\mod{p_1p_2})$, where $h$ is the class number of $\Q(\sqrt{-p_1p_2})$. It is easy to see by induction that 
$\ord_{N}(p)=\phi(N)/2$ for all pairs $(p,N)$ in Table~\ref{Subsec1Tab2}. 
Furthermore, since $(p^{p_1^{m-1}p_2^{n-1}})^{\frac{p_1-1}{2}\frac{p_2-1}{2}}\equiv p^{\frac{p_1-1}{2}\frac{p_2-1}{2}}\,(\mod{p_1p_2})$, 
the condition $bp^{\frac{f-h}{2}}\equiv 2\,(\mod{p_1p_2})$ can be rewritten as 
$bp^{\frac{p_1-1}{2}\frac{p_2-1}{2}}\equiv 2p^{\frac{h}{2}}\,(\mod{p_1p_2})$, 
which is independent of  $m$ and $n$. There are only these three series satisfying the 
conditions of Corollary~\ref{Subsec1Cor1} when $p_1\le 10^{7}$. 
\begin{table}[h]
\begin{center}
\caption{\label{Subsec1Tab2}
Generalizations of the strongly regular graphs in Table 1.  The parameters $\lambda$ and $\mu$ of the strongly regular graphs can be computed by $\lambda=s+r+sr+k$ and $\mu=k+sr$. }
$$
\begin{array}{|c||c|c|c|c|c|c|c|c|}
\hline
\mbox{No.} &\mbox{$p$} & \mbox{$N$}&\mbox{$h$}& \mbox{$b$} & \mbox{$v$}&\mbox{$k$}
&\mbox{$r$, $s$}\\
\hline \hline 
\mbox{7} &\mbox{$2$} & \mbox{$3^m\cdot 5^n$}&\mbox{$2$}& \mbox{$1$} & \mbox{$2^{4\cdot 3^{m-1}\cdot 5^{n-1}}$} & \mbox{$\frac{2^{4\cdot 3^{m-1}\cdot 5^{n-1}}-1}{15}$}
&\mbox{$r=\frac{8\cdot 2^{2\cdot 3^{m-1}\cdot 5^{n-1}-1}-1}{15}$}\\
\mbox{}&\mbox{}&\mbox{}&\mbox{}&\mbox{}&\mbox{}&\mbox{}
&\mbox{$s=\frac{-7\cdot 2^{2\cdot 3^{m-1}\cdot 5^{n-1}-1}-1}{15}$}\\
\hline
\mbox{8} &\mbox{$3$} & \mbox{$5^m\cdot 7^n$}&\mbox{$2$}& \mbox{$-1$} & \mbox{$3^{12\cdot 5^{m-1}\cdot7^{n-1}}$} & \mbox{$\frac{3^{12\cdot 5^{m-1}\cdot7^{n-1}}-1}{35}$}
&\mbox{$r=\frac{17\cdot 3^{6\cdot 5^{m-1}\cdot7^{n-1}-1}-1}{35}$}\\
\mbox{}&\mbox{}&\mbox{}&\mbox{}&\mbox{}&\mbox{}&\mbox{}
&\mbox{$s=\frac{-18\cdot 3^{6\cdot 5^{m-1}\cdot7^{n-1}-1}-1}{35}$}\\
\hline
\mbox{9} &\mbox{$3$}&\mbox{$17^m\cdot 19^n$}&\mbox{$4$}&\mbox{$-1$}&\mbox{$3^{144\cdot 17^{m-1}\cdot 19^{n-1}}$} & \mbox{$\frac{3^{144\cdot 17^{m-1}\cdot 19^{n-1}}-1}{323}$}&\mbox{$r=\frac{161\cdot 3^{72\cdot 17^{m-1}\cdot 19^{n-1}-2}-1}{323}$}\\
\mbox{}&\mbox{}&\mbox{}&\mbox{}&\mbox{}&\mbox{}&\mbox{}
&\mbox{$s=\frac{-162\cdot 3^{72\cdot 17^{m-1}\cdot 19^{n-1}-2}-1}{323}$}\\
\hline
\end{array}
$$
%\hspace{0.8cm}
\end{center}
\end{table}
\end{example}
%%%%%%%%%%%%%%%%%%%%%%%%%%%%%%%%%%%%%%%%%%%%%%%%%%%%%%%%%
%%%%%%%%%%%%%%%%%%%%%%%%%%%%%%%%%%%%%%%%%%%%%%%%%%%%%%%%%
\subsection{Skew Hadamard difference sets from unions of cyclotomic classes of order $N=2p_1^m$}
In this subsection, we assume that 

\begin{enumerate}

\item  $p_1\equiv 3\pmod{8}$, ($p_1\neq 3$),

\item  $N=2p_1^m$,

\item  $1+p_1=4p^h$, where $h$ is the class number of $\Q(\sqrt{-p_1})$,

\item  $p$ is a prime such that $\ord_{N}(p)=\phi(p_1^m)/2$. 

\end{enumerate}

Let $q=p^f$, where $f=\ord_N(p)$. 
Let $\zeta_{q-1}=\exp(\frac{2\pi i}{q-1})$ and $\mathfrak{P}$ be a 
prime ideal in $\Z[\zeta_{q-1}]$ lying over $p$. Then,  
$\Z[\zeta_{q-1}]/\mathfrak{P}$ is the finite field of order $q$ and 
written as $
\Z[\zeta_{q-1}]/\mathfrak{P}=\{\overline{\zeta}_{q-1}^i\,|\,0\le i\le q-2\}\cup \{\overline{0}\}$, where 
$\overline{\zeta}_{q-1}=\zeta_{q-1}+\mathfrak{P}$. 
Hence, $\gamma:=\overline{\zeta}_{q-1}$ is a primitive element of 
$\F_q=\Z[\zeta_{q-1}]/\mathfrak{P}$. 
Let $\omega_{\mathfrak{P}}$ be the 
Teichm\"{u}ller character  of $\F_q$. Then, 
$\omega_{\mathfrak{P}}(\gamma)=\zeta_{q-1}$. 
Put $\chi_N:=\omega_{\mathfrak{P}}^\frac{q-1}{N}$. Then 
$\chi_N$ is a multiplicative character of order $N$ of $\F_q$. For this 
$\chi_N$, by the results of \cite{L97}, we have
\begin{equation}\label{Subsec2Eq1}
G(\chi_{N}^2)=G(\chi_{p_1^m})=p^{\frac{f-h}{2}}\left(\frac{b+c\sqrt{-p_1}}{2}\right), 
\end{equation}
where $b,c\not \equiv 0\,(\mod{p})$, $b^2+c^2p_1=4p^h$, and 
$bp^\frac{f-h}{2}\equiv -2\,(\mod{p_1})$. By our assumption that 
$1+p_1=4p^h$, we have $b,c\in \{-1,1\}$, where the sign of $c$ depends on the choice of $\mathfrak{P}$. In particular, 
in \cite{FX113}, it was shown that $bc\equiv -\sqrt{-p_1}\,(\mod{\mathfrak{P}})$. 
On the other hand, since $1+p_1=4p^h$, we have 
$(1+\sqrt{-p_1})(1-\sqrt{-p_1})\in \mathfrak{P}$, from which it follows that $1+\sqrt{-p_1}\in \mathfrak{P}$ or 
$1-\sqrt{-p_1}\in \mathfrak{P}$ for any prime ideal 
$\mathfrak{P}$ in $\Q(\zeta_{q-1})$ lying over $p$. We may choose a prime ideal $\mathfrak{P}$ such that $1+\sqrt{-p_1} \in  \mathfrak{P}$. 
Then, $bc\equiv -\sqrt{-p_1}\,(\mod{\mathfrak{P}})$ with $b,c\in \{-1,1\}$ implies that $bc=1$. 
From now on, we fix this choice of $\mathfrak{P}$. 

Let $C_i=\gamma^i \langle \gamma^N\rangle$, where $0\leq i\leq N-1$, and $\gamma$ is the 
fixed primitive element of $\F_q$ as above. It is clear that $-1\in C_0$ or 
$-1\in C_{p_1^m}$ depending on whether $p\equiv 1\,(\mod{4})$ or 
$p\equiv 3\,(\mod{4})$. 

Let 
\[
J=\langle p\rangle \cup 2\langle p\rangle
\cup \{0\}\, (\mod{2p_1})
\]
and 
define 
\[
D=\bigcup_{i=0}^{p_1^{m-1}-1}\bigcup_{j\in J}C_{2i+p_1^{m-1}j}. 
\]  
From the facts that $2$ is a nonsquare of $\F_{p_1}$ and that the reduction of $
\langle p\rangle \le \Z_{N}^\ast$ modulo $2p_1$ 
is the subgroup of index $2$ of $\Z_{2p_1}^\ast$ we deduce that 
$J\,(\mod{p_1})=\Z_{p_1}$, and $D=-D$ or $D\cap -D=\emptyset$ 
depending on whether $p\equiv 1\,(\mod{4})$ or 
$p\equiv 3\,(\mod{4})$.
\begin{theorem}\label{Subsec2Thm1}
The size of the set $\{\psi(\gamma^a D)\,|\,a=0,1,\ldots,q-2\}$ is precisely two.  
\end{theorem}
\proof 
Set $A:=(-1)^{\frac{p-1}{2}(m-1)}p^{\frac{f-1}{2}-h}\sqrt{p^\ast}$
and $B:=(-1)^{\frac{p-1}{2}\frac{f-1}{2}}p^{\frac{f-1}{2}}\sqrt{p^\ast}$, where 
$\sqrt{p^\ast}=\sqrt{(-1)^{\frac{p-1}{2}}p}$.

First of all, we note that $(-1)^{\frac{f-1}{2}}=(-1)^{m-1}$ since 
$p_1\equiv 3\,(\mod{8})$.  It follows that $p^hA=B$.

Secondly, since $D$ is a union of cyclotomic classes of order $N$, we have $\{\psi(\gamma^a D)\,|\,a=0,1,\ldots,q-2\}=\{\psi(\gamma^a D)\,|\,a=0,1,\ldots, N-1\}$.

It is sufficient to evaluate the sums   
\[
T_a:=N\cdot \psi(\gamma^a D)=
\sum_{\ell=0}^{2p_1^m-1}G(\chi_{2p_{1}^m}^{\ell}) \sum_{i=0}^{p_1^{m-1}-1}
\sum_{j\in J}\chi_{2p_{1}^m}^{-\ell}(\gamma^{a+2i+p_1^{m-1}j}), 
\]
 where $a=0,1,\ldots, N-1$.

For $\ell=0$, by noting that $G(\chi_{2p_{1}^m}^{0})=-1$, we have 
\[
G(\chi_{2p_{1}^m}^{0}) \sum_{i=0}^{p_1^{m-1}-1}
\sum_{j\in J}\chi_{2p_{1}^m}^{0}(\gamma^{a+2i+p_1^{m-1}j})=-p_1^{m}. 
\]

For $\ell=2 h$ but $h\not\equiv 0\,(\mod{p_1})$,  since  $J\,(\mod{p_1})=\Z_{p_1}$, 
we have 
\[
\sum_{j\in J}
\chi_{2p_{1}^m}^{-2h}(\gamma^{p_1^{m-1}j})=\sum_{j\in J}
\chi_{p_{1}}^{-h}(\gamma^{j})=0.
\]

For $\ell=p_1 h$ but  $h\not\equiv 0\,(\mod{p_1^{m-1}})$,  
we have 
\[
\sum_{i=0}^{p_1^{m-1}-1}
\chi_{2p_{1}^m}^{-p_1 h}(\gamma^{2i})=
\sum_{i=0}^{p_1^{m-1}-1}
\chi_{p_{1}^{m-1}}^{-h}(\gamma^{i})=0. 
\]

For $\ell=p_1^{m}$,  since $G(\chi_{2p_1^m}^{p_1^m})=B$ by Theorem~\ref{Sec2Thm2}, 
we have 
\[
G(\chi_{2p_1^m}^{p_1^m})\sum_{i=0}^{p_1^{m-1}-1}\sum_{j\in J}
\chi_{2p_{1}^m}^{p_1^m }(\gamma^{a+2i+p_1^{m-1}j})=Bp_1^{m-1}(-1)^a. 
\]

For the remaining cases, we evaluate the sum  
\[
\sum_{\ell\in \langle p\rangle}
G(\chi_{2p_{1}^m}^{\ell}) \sum_{i=0}^{p_1^{m-1}-1}\sum_{j\in J}
\chi_{2p_{1}^m}^{-\ell}(\gamma^{a+2i+p_1^{m-1}j})
+\sum_{\ell\in - \langle p\rangle}
G(\chi_{2p_{1}^m}^{\ell}) \sum_{i=0}^{p_1^{m-1}-1}\sum_{j\in J}
\chi_{2p_{1}^m}^{-\ell}(\gamma^{a+2i+p_1^{m-1}j}). 
\]
By Theorem~\ref{Sec2Thm2}, we have 
\[
G(\chi_{2p_1^m}^\ell)=A\left(\frac{b+c\sqrt{-p_1}}{2}\right)^2
\] 
for $\ell\in \langle p\rangle$, where $b,c$ are the same as in the 
evaluation (\ref{Subsec2Eq1}) of $G(\chi_{p_1^m})$. By the choice of $\mathfrak{P}$, 
it is expanded  as 
\[
G(\chi_{2p_1^m}^\ell)=A\left(\frac{1-p_1+2\sqrt{-p_1}}{4}\right). 
\]

Since $\chi_{2p_1^m}^\ell (-1)\overline{\sqrt{p^\ast}}=\sqrt{p^\ast}$ for any odd 
$\ell$ by the assumption $p_1\equiv 3\,(\mod{8})$, i.e., $f$ is odd, 
the above sum is reformulated as 
\begin{eqnarray*}
& &A\left(\frac{1-p_1+2\sqrt{-p_1}}{4}\right)\sum_{\ell\in \langle p\rangle}
 \sum_{i=0}^{p_1^{m-1}-1}\sum_{j\in J}\chi_{2p_{1}^m}^{-\ell}(\gamma^{a+2i+p_1^{m-1}j})\\
& &\hspace{2cm}+A\left(\frac{1-p_1-2\sqrt{-p_1}}{4}\right)\sum_{\ell\in -\langle p\rangle} \sum_{i=0}^{p_1^{m-1}-1}\sum_{j\in J}
\chi_{2p_{1}^m}^{-\ell}(\gamma^{a+2i+p_1^{m-1}j}). 
\end{eqnarray*}

Note that $\langle p\rangle $ can be written as 
$\{x+2p_1y\,|\,x\in \langle p\rangle (\mod{2p_1}),y\in \{0,1,\ldots,p_1^{m-1}-1\}\}$. 
Furthermore, there is a unique $i\in \{0,1,\ldots ,p_1^{m-1}-1\}$ such that $a+2i\equiv 0\,(\mod{p_1^{m-1}})$; we write 
$a+2i=p_1^{m-1}i_a$. 
Then, the above sum is rewritten as 
\begin{eqnarray*}
& &A\left(\frac{1-p_1+2\sqrt{-p_1}}{4}\right)\sum_{x\in\langle p\rangle (\mod{2p_1})}\sum_{y=0}^{p_1^{m-1}-1}
 \sum_{i=0}^{p_1^{m-1}-1}\sum_{j\in J}\chi_{2p_{1}^m}^{-x}(\gamma^{a+2i+p_1^{m-1}j})\chi_{2p_{1}^m}^{-2p_1y}(\gamma^{a+2i+p_1^{m-1}j})\\
& &\hspace{0.5cm}+A\left(\frac{1-p_1-2\sqrt{-p_1}}{4}\right)\sum_{x\in -\langle p\rangle (\mod{2p_1})}\sum_{y=0}^{p_1^{m-1}-1} \sum_{i=0}^{p_1^{m-1}-1}\sum_{j\in J}
\chi_{2p_{1}^m}^{-x}(\gamma^{a+2i+p_1^{m-1}j})\chi_{2p_{1}^m}^{-2p_1y}(\gamma^{a+2i+p_1^{m-1}j})\\
&=&p_1^{m-1}A\left(\frac{1-p_1+2\sqrt{-p_1}}{4}\right)\sum_{x\in\langle p\rangle (\mod{2p_1})}
\sum_{j\in J}\chi_{2p_{1}}^{-x}(\gamma^{i_a+j})\\
& &\hspace{1.5cm}+p^{m-1}A\left(\frac{1-p_1-2\sqrt{-p_1}}{4}\right)\sum_{x\in -\langle p\rangle (\mod{2p_1})}\sum_{j\in J}
\chi_{2p_{1}}^{-x}(\gamma^{i_a+j})\\
\\
&=&p_1^{m-1}A\left(\frac{1-p_1+2\sqrt{-p_1}}{4}\right)\left(\sum_{x\in\langle p\rangle (\mod{2p_1})}\chi_{2p_{1}}^{-x}(\gamma^{i_a})\right)\left(
\sum_{j\in J}\chi_{2p_{1}}^{-j}(\gamma)\right)\\
& &\hspace{1.5cm}+p^{m-1}A\left(\frac{1-p_1-2\sqrt{-p_1}}{4}\right)\left(\sum_{x\in \langle p\rangle (\mod{2p_1})}\chi_{2p_{1}}^{x}(\gamma^{i_a})\right)\left(\sum_{j\in J}
\chi_{2p_{1}}^{j}(\gamma)\right).
\end{eqnarray*}
Put
$X_a:=\sum_{j\in J}\chi_{2p_{1}}^{-j}(\gamma)$ and 
 $Y_a:=\sum_{x\in\langle p\rangle (\mod{2p_1})}\chi_{2p_{1}}^{-x}(\gamma^{i_a})$. Let $\eta$ be the quadratic character of $\F_{p_1}$ and 
$\psi_{p_1}$ be the canonical additive character of $\F_{p_1}$. 
Noting that $2$ is a nonsquare in $\F_{p_1}$.  For
$i\in \Z_{2p_1}\setminus\{0,p_1\}$ it holds that 
\begin{eqnarray*}
\sum_{x\in\langle p\rangle (\mod{2p_1})}\chi_{2p_{1}}^{-x}(\gamma^i)&=&
\sum_{x\in\langle p\rangle (\mod{2p_1})}\chi_2^{-ix}(\gamma)\chi_{p_{1}}^{\frac{p_1-1}{2}ix}(\gamma)\\
&=&
(-1)^i\frac{1}{2}\sum_{x\in \F_{p_1}^\ast}(1+\eta(x))\psi_{p_1}(-2^{-1}ix)\\
&=&
(-1)^i\frac{-1+\eta(-2^{-1}i)G(\eta)}{2}=(-1)^i
\frac{-1+\eta(i)\sqrt{-p_1}}{2}.  
\end{eqnarray*}
Hence, we have 
\begin{eqnarray*}
X_a&=&\sum_{j\in \langle p\rangle\,(\mod{2p_1})}\chi_{2p_{1}}^{-j}(\gamma)+
\sum_{j\in 2\langle p\rangle\,(\mod{2p_1})}\chi_{2p_{1}}^{-j}(\gamma)+1\\
&=&
\frac{1-\sqrt{-p_1}}{2}+
\frac{-1-\sqrt{-p_1}}{2}+1=1-\sqrt{-p_1}
\end{eqnarray*}
and 
\begin{eqnarray*}
Y_a=(-1)^{i_a}
\frac{-1+\eta(i_a)\sqrt{-p_1}}{2}, \; i_a\neq 0, p_1. 
\end{eqnarray*}
Thus, we obtain
\begin{eqnarray*}
& &T_a+p_1^{m}\\
&=&Bp_1^{m-1}(-1)^a+\frac{p_1^{m-1}A}{4}\left((1-p_1+2\sqrt{-p_1})(1-\sqrt{-p_1})Y_a
+(1-p_1-2\sqrt{-p_1})(1+\sqrt{-p_1})\overline{Y_a}\right). 
\end{eqnarray*}

We compute $T_a+p_1^m$ by considering the following six cases:  
\begin{itemize}
\item[(i)] $i_a=0$: 
In this case, we have $a\equiv 0\,(\mod{2})$, $Y_a=\frac{p_1-1}{2}$, and 
$
T_a+p_1^{m}=p_1^{m-1}(\frac{A}{4}(p_1^2-1)+B). 
$
\item[(ii)] $i_a=p_1 $: In this case, we have  
$a\equiv 1\,(\mod{2})$, $Y_a=-\frac{p_1-1}{2}$, and 
$
T_a+p_1^{m}=p_1^{m-1}(\frac{A}{4}(-p_1^2+1)-B).
$
\item[(iii)] $i_a\in \langle p\rangle $: In this case, we have  $a\equiv 1\,(\mod{2})$,  $Y_a=\frac{1-\sqrt{-p_1}}{2}$,  and 
$
T_a+p_1^{m}=p_1^{m-1}(\frac{A}{4}(p_1^2+2p_1+1)-B). 
$
\item[(iv)] $i_a\in -\langle p\rangle $: In this case, we have 
$a\equiv 1\,(\mod{2})$,  $Y_a=\frac{1+\sqrt{-p_1}}{2}$, and 
$
T_a+p_1^{m}=p_1^{m-1}(\frac{A}{4}(1-p_1^2)-B).
$
\item[(v)] $i_a\in 2\langle p\rangle $: In this case, we have 
$a\equiv 0\,(\mod{2})$, $Y_a=-\frac{1+\sqrt{-p_1}}{2}$, and 
$
T_a+p_1^{m}=p_1^{m-1}(\frac{A}{4}(p_1^2-1)+B). 
$
\item[(vi)] $i_a\in -2\langle p\rangle $: In this case, we have $a\equiv 0\,(\mod{2})$, $Y_a=\frac{-1+\sqrt{-p_1}}{2}$, and 
$
T_a+p_1^{m}=p_1^{m-1}(\frac{A}{4}(-p_1^2-2p_1-1)+B).
$
\end{itemize}
By the  assumption that 
$1+p_1=4p^h$ and the fact that $p^hA=B$, it is easily checked that 
$T_a+p_1^{m}$, $a=0,1,\ldots, N-1$, take precisely two values. The proof is now complete. 
\qed

\begin{corollary}\label{Subsec2Cor1}
The set $D$ is a skew Hadamard difference set or 
a Paley type partial difference set according as $p\equiv 3\,(\mod{4})$ or $p\equiv 1\,(\mod{4})$. 
\end{corollary}
\proof
By Theorem~\ref{Subsec2Thm1}, the set $\{\psi(\gamma^a D)\,|\,a=0,1,\ldots,q-2\}$ has precisely two elements, which are 
\[
\frac{1}{N}\left(-p_1^{m}\pm p_1^{m-1}\left(\frac{A}{4}(p_1^2-1)+B\right)\right)
=\frac{1}{2}\left(-1\pm (-1)^{\frac{(p-1)(m-1)}{2}}\sqrt{(-1)^{\frac{p-1}{2}}p^f}\right).
\]
By Lemma~\ref{Sec3Le1}, the assertion of the corollary follows immediately. 
\qed

The construction of skew Hadamard difference sets and 
Paley type partial difference sets  given in this subsection is a generalization of Theorem~\ref{IntroThm3} (ii) \cite{FX113}. 
In particular, in \cite{FX113},  one example of skew Hadamard difference sets with parameters 
$(p,N,h,b,v)=(3,2\cdot 11,1,1,3^{5})$ was given. Unfortunately, 
we can not generalize this example to $N=2\cdot 11^m$ because $p=3$ does not satisfy the condition $\ord_{N}(p)=\phi(N)/2$ for $m>1$. Below are some infinite series of skew Hadamard difference sets and Paley type partial difference 
sets obtained by Corollary~\ref{Subsec2Cor1}. 
\begin{example}
Table~\ref{Subsec2Tab1} shows all possible 
skew Hadamard difference sets and Paley type partial difference sets 
obtained by applying Corollary~\ref{Subsec2Cor1} to all  $p_1\le 10^6$ except for
the case when $p_1=11$ and $m=1$.
In particular, the $3$rd case of Table~\ref{Subsec2Tab1} gives skew Hadamard difference sets and 
the other cases give Paley type partial difference 
sets. Here, 
the parameters $p,N,h,b$ satisfy the conditions of Corollary~\ref{Subsec2Cor1}, i.e., $p$ is a prime such that 
$[\Z_{N}^\ast:\langle p\rangle]=2$, $1+p_1=4p^h$, and 
$bp^{\frac{f-h}{2}}\equiv -2\,(\mod{p_1})$, where $h$ is the class number of $\Q(\sqrt{-p_1})$. Note that it is easy to prove by induction that 
$\ord_{N}(p)=\phi(N)/2$ for all pairs $(p,N)$ in Table~\ref{Subsec2Tab1}. 
Furthermore, since 
\[
p^{\frac{p_1^{m-1}(p_1-1)/2-h}{2}}\equiv 
p^{\frac{p_1^{m-1}(p_1+1)}{4}-\frac{p_1^{m-1}-1}{2}-\frac{h+1}{2}}\, (\mod{p_1}), 
\]
the condition $bp^{\frac{f-h}{2}}\equiv -2\,(\mod{p_1})$ can be rewritten as 
$
bp^{\frac{p_1-1-2h}{4}}\equiv -2\,(\mod{p_1}), 
$
which is independent of $m$. 
\begin{table}[h]
\begin{center}
\caption{\label{Subsec2Tab1}
Some Paley type partial difference sets and skew Hadamard difference sets  obtained by Corollary~\ref{Subsec2Cor1}. }
$$
\begin{array}{|c||c|c|c|c|c|c|c|}
\hline
\mbox{No.} &\mbox{$p$} & \mbox{$N$}&\mbox{$h$}& \mbox{$b$} & \mbox{$v$}\\
\hline \hline 
\mbox{1} &\mbox{$5$} & \mbox{$2\cdot 19^m$}&\mbox{$1$}& \mbox{$1$} & \mbox{$5^{9\cdot 19^{m-1}}$} \\
\hline
\mbox{2} &\mbox{$17$} & \mbox{$2\cdot 67^m$}&\mbox{$1$}& \mbox{$1$} & \mbox{$17^{33\cdot 67^{m-1}}$}\\
\hline
\mbox{3} &\mbox{$3$}&\mbox{$2\cdot 107^m$}& \mbox{$3$}&\mbox{$1$}& \mbox{$3^{53\cdot 107^{m-1}}$} \\
\hline
\mbox{4} &\mbox{$41$}&\mbox{$2\cdot 163^m$}& \mbox{$1$}&\mbox{$1$}&\mbox{$41^{81\cdot 163^{m-1}}$} \\
\hline
\mbox{5} &\mbox{$5$}&\mbox{$2\cdot 499^m$}& \mbox{$3$}&\mbox{$1$}&\mbox{$5^{249\cdot 499^{m-1}}$} \\
\hline
\end{array}
$$
%\hspace{0.8cm}
\end{center}
\end{table}
\end{example}
\section{Concluding remarks}
In this paper, we have given two constructions of strongly regular graphs and skew Hadamard difference sets, 
which are generalizations of those given by the first and third authors \cite{FX111,FX113}. As a consequence, we obtain three infinite series of strongly regular graphs with new parameters and a family of skew Hadamard difference sets in $(\F_q,+)$, where $q=3^{53\cdot 107^{m-1}}$. The results on strongly regular graphs have implications on association schemes.  

Given a $d$-class (symmetric) association scheme $(X, \{R_\ell\}_{0\leq
\ell\leq d})$, we can take the union of classes to form graphs with
larger edge sets (this process is called a {\em{fusion}}), but it is
not necessarily guaranteed that the fused collection of graphs will
form an association scheme on $X$. If an association scheme has the
property that any of its fusions is also an association scheme, then
we call the association scheme {\em{amorphic}}. A well-known and
important example of amorphic association schemes is given by the
cyclotomic association schemes on $\F_q$ when the cyclotomy is uniform \cite{BMW82}.

In \cite{IP94}, A.V. Ivanov conjectured that if each nontrivial relation in an 
association scheme is strongly regular, then the association scheme 
must be amorphic. This conjecture turned out to be false. A first counterexample was found by 
Van Dam \cite{D00} in the case when the association scheme is imprimitive. 
Afterwards, Van Dam \cite{D03} and Ikuta and Munemasa \cite{IM10} gave more counterexamples in the case when 
the association scheme is primitive. However, there had been known only a few counterexamples in the primitive case. 
Recently, in \cite{FX112}, the authors generalized the counterexamples of Van Dam and 
Ikuta-Munemasa into infinite series using strongly regular Cayley graphs based on index $2$ Gauss sums 
of type $N=p_1^m$ and type $N=p_1^mp_2$. Our generalization (Corollary~\ref{Subsec1Cor1}) of the second construction in \cite{FX111} 
produces further new counterexamples to Ivanov's conjecture and association schemes with very interesting properties.
More precisely, under the same assumptions as in Corollary~\ref{Subsec1Cor1},
define 
\[
D_k=\bigcup_{i=0}^{p_1^{m-1}-1}\bigcup_{j=0}^{p_2^{n-1}-1}
C_{p_2^n i+p_1^m j+p_1^{m-1}p_2^{n-1}k}
\]
for each $0\le k\le p_1p_2-1$. 
Let $R_0=\{(x,x)\,|\,x\in \F_q\}$ and 
\[
R_k:=\{(x,y)\,|\,x,y\in \F_q,x-y\in D_{k-1}\}. 
\] 
Then, one can similarly prove that $(\F_q,\{R_k\}_{0\le k\le p_1p_2})$ is a pseudocyclic and non-amorphic  
association scheme in which every nontrivial relation is a strongly regular graph. Table~\ref{Subsec1Tab2} yields 
three new infinite series of pseudocyclic and non-amorphic association schemes, where each of the nontrivial relations is strongly regular.  
Moreover, further fusion schemes of these association schemes are possible by applying Corollary 3.2 and Theorem 4.1 of \cite{IM10}. 
In particular, Examples 1 and 2 of \cite{IM10} are generalized into an infinite series by using the above association scheme with 
$p=2$, $b=1$, $(p_1,p_2)=(5,3)$, and $h=2$. 

\section*{Acknowledgements} 
The work of Tao Feng was supported in part by the Fundamental Research Funds for the central universities. The work of K. Momihara was supported by JSPS under Grant-in-Aid for Research Activity Start-up 23840032. The work of Qing Xiang was supported in part by NSF Grant DMS 1001557 and by the Y. C. Tang disciplinary development fund of Zhejiang University.

\end{document}